\documentclass{amsart}

\usepackage{amsmath,amssymb,amsthm}

\usepackage{graphicx,tikz, caption}

\newtheorem*{thm}{Theorem}

\newtheorem{corollary}{Corollary}

\begin{document}

\title[]{An Almgren monotonicity formula \\for discrete harmonic functions}

\author[]{Mariana Smit Vega Garcia}
\address{Department of Mathematics, Western Washington University, Bellingham, WA 98225}
\email{smitvem@wwu.edu}

\author[]{Stefan Steinerberger}
\address{Department of Mathematics, University of Washington, Seattle, WA 98195, USA} 
\email{steinerb@uw.edu}

\date{}
\thanks{The authors have been supported by NSF DMS-2054282 (MSVG) and DMS-2123224 (SS)}

\begin{abstract} The celebrated Almgren monotonicity formula for harmonic functions $u:\mathbb{R}^n \rightarrow \mathbb{R}$ says that its $L^2$-energy concentrated on a sphere of radius $r$, when measured in a suitable sense, is non-decreasing: if $u$ oscillates at a certain scale, it has even larger oscillations at a larger scale.
We prove a discrete analogue of the Almgren monotonicity formula for harmonic functions on infinite combinatorial graphs $G=(V,E)$. Some applications are discussed.
\end{abstract}

\maketitle

\section{Introduction and Results}
\subsection{Introduction}
We will show that the classical Almgren monotonicity formula for harmonic functions $u:\mathbb{R}^n \rightarrow \mathbb{R}$ has an extension to the setting of discrete harmonic functions on combinatorial graphs. The celebrated Almgren monotonicity formula \cite{almgren, almgren2} is a cornerstone in the study of harmonic functions. It also plays a crucial role in studying unique continuation, and has been used extensively in free boundary
problems \cite{ACS, PSU, sp}. It is usually stated as follows: let $u:\mathbb{R}^n \rightarrow \mathbb{R}$ be a harmonic function and let $B_r$ denote a ball of radius $r$, then
$$ N(r) = \frac{r \int_{B_r} |\nabla u|^2 dx}{\int_{\partial B_r} u(x)^2 d\sigma} \qquad \mbox{is non-decreasing}.$$
One often lets $r \rightarrow 0$ to deduce information from the limit. Conversely, the presence of oscillation implies the existence of larger oscillations at a larger scale.
Using integration by parts, the Dirichlet energy can be written as
\begin{align*} \int_{B_r} |\nabla u|^2 dx & = \int_{\partial B_r} u \frac{\partial u}{\partial n} d\sigma \\
&= \frac{1}{2} \frac{d}{dr} \int_{\partial B_r} u(x)^2 d\sigma-\frac{n-1}{r}\int_{\partial B_r}u(x)^2d\sigma.
\end{align*}
Therefore, Almgren's monotonicity can be equivalently written as saying that
$$ N(r) = r \frac{d}{dr} \log\left(\int_{\partial B_r} u(x)^2 d\sigma \right) \quad \mbox{is non-decreasing.}$$
This implies that if a harmonic function is large in $L^2$ on a sphere, then it will be even larger on spheres with the same center and larger radius.
After suitably adapting the functional to the setting of combinatorial graphs, we will obtain an analogous identity with the same consequences.

\subsection{Discrete Almgren Monotonicity.}
We can now introduce a discrete analogue.  Recall that a function $u:V \rightarrow \mathbb{R}$ is harmonic if, for each vertex $x \in V$,
$$ \sum_{(x, y) \in E} (u(x) - u(y)) = 0 \qquad \mbox{or, equivalently} \qquad  u(x) = \frac{1}{\deg(x)} \sum_{(x,y) \in E} u(y).$$
Such functions have been intensively studied and have a rich theory, see, for example \cite{graph0, graph1, graph2, graph3}.
We assume we are given a connected graph $G=(V,E)$. The graph may be finite (in which case our result will be trivial: harmonic functions are constant). 
If the graph is infinite, we require that it is locally finite: each vertex has finite degree. We interpret Almgren's functional quite literally and regard
$$ \int_{\partial B_r} u(x)^2 dx \qquad \mbox{as} \qquad \sum_{d(x,y)=k} u(y)^2,$$
where $x \in V$ is an arbitrary vertex.
The differentiation in $r$ is replaced by a discrete difference between $k+1$ and $k$. There is an important conceptual change: vertices with a higher degree have a bigger impact on their neighborhood and this has to be accounted for. Fixing the base vertex $x \in V$ induces a partition of the vertex set by distance from $x$ (see Fig. 1). We introduce the in-degree of a vertex $v \in V$ 
$$ d_{\mbox{\tiny in}}(v) = \# \left\{w \in V: (w,v) \in E \quad \mbox{and} \quad d(w,x) = d(v,x) -1 \right\}$$
and, completely analogously, the out-degree as
$$ d_{\mbox{\tiny out}}(v) = \# \left\{w \in V: (w,v) \in E \quad \mbox{and} \quad d(w,x) = d(v,x) +1 \right\}.$$
For any $v \notin x$, we have $d_{\mbox{\tiny in}}(v) \geq 1$ as well as $d_{\mbox{\tiny in}}(v) + d_{\mbox{\tiny out}}(v) \leq \deg(v)$.
\begin{center}
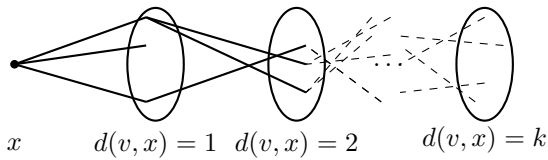
\begin{figure}[h!]
\begin{tikzpicture}[scale=1.25]
\filldraw (-1,0) circle (0.04cm);
\draw[thick] (0.5,0) ellipse (0.3cm and 0.6cm);
\draw[thick] (2,0) ellipse (0.3cm and 0.6cm);
\node at (3,0) {$\dots$};
\draw[thick] (4,0) ellipse (0.3cm and 0.6cm);
\draw [thick] (-1,0) -- (0.9-0.5, 0.2);
\draw [thick] (-1,0) -- (0.9-0.5, -0.4);
\draw [thick] (-1,0) -- (0.9-0.5, 0.5);
\draw [thick] (0.9-0.5, 0.5) -- (2.1, 0);
\draw [thick] (0.9-0.5, 0.5) -- (2.1, -0.3);
\draw [thick] (0.9-0.5, -0.4) -- (2.1, 0.2);
\draw [dashed] (3, 0.5) -- (2.1, 0);
\draw [dashed] (2.8, 0.5) -- (2.1, -0.3);
\draw [dashed] (2.9, -0.4) -- (2.1, 0.2);
\draw [dashed] (3, 0.1) -- (2.1, 0);
\draw [dashed] (2.8, 0.3) -- (2.1, -0.3);
\draw [dashed] (2.9, -0.4) -- (2.1, 0.2);
\node at (-1, -0.85) {$x$};
\node at (0.5, -0.85) {$d(v,x) = 1$};
\node at (2, -0.85) {$d(v,x) = 2$};
\node at (4, -0.8) {$d(v,x) = k$};
\draw [dashed] (4, 0.5) -- (3.1, 0);
\draw [dashed] (4, -0.2) -- (3.1, -0.3);
\draw [dashed] (4.2, 0.2) -- (3.1, 0.3);
\draw [dashed] (3.9, -0.4) -- (3.1, 0.2);
\end{tikzpicture}
\caption{Ordering vertices by distance from a fixed vertex.}
\end{figure}
\end{center}

\begin{thm}[Discrete Almgren Monotonicity Formula]\label{T:Almgren} Let $G =(V,E)$ be locally finite, let $x \in V$ and let $u:V \rightarrow \mathbb{R}$ be harmonic. Then, for $k \geq 0$, 
$$N(k) =  \sum_{d(x,y)=k+1} d_{\emph{\tiny in}}(y) \cdot u(y)^2  - \sum_{d(x,y)=k} d_{\emph{\tiny out}}(y) \cdot u(y)^2 $$
is non-negative and satisfies $N(k+1) \geq N(k).$
\end{thm}

The result has an immediate extension to weighted graphs, see \S 2.5.
If there is a dramatic change in the (suitably normalized) $\ell^2$-energy on a sphere, then this is indicative of even larger fluctuations at larger spheres.  When $u \equiv 1$ is constant, we have $N \equiv 0$ since the number of incoming edges in $\left\{y: d(x,y)=k+1 \right\}$ is exactly the same as the number of outgoing edges from $\left\{y: d(x,y)=k \right\}$. The inequality is sharp in the sense that $N(k+1)-N(k)$ may tend to 0, see \S 2.4.

\subsection{Doubling estimates.} The continuous Almgren monotonicity formula immediately leads to doubling estimates. Something similar is true in the discrete setting. Let $0 \leq a \leq b < \infty$ be two integers, then the quantities $N(a)$ and $N(b)$ are related to the growth of the harmonic function.
We say that a graph is \textit{locally expansive} in $\left\{v \in V: a \leq d(x,v) \leq b\right\}$ if for all vertices satisfying $a \leq d(v,x) \leq b$, we have
$$ d_{\mbox{\tiny out}}(v) \geq d_{\mbox{\tiny in}}(v).$$
Likewise, we say that the graph is \textit{locally contractive} if $ d_{\mbox{\tiny out}}(v) \leq d_{\mbox{\tiny in}}(v)$ for every vertex in that region. In either of these cases, we can obtain bounds on the growth in terms of $N(k)$ via a simple telescoping argument.

\begin{corollary}
Let $G =(V,E)$ be locally finite, let $x \in V$ and let $u:V \rightarrow \mathbb{R}$ be a harmonic function. If $G$ is locally expansive in $\left\{v \in V: a \leq d(x,v) \leq b\right\}$, then
$$ \sum_{d(x,y)=b+1} d_{\emph{\tiny in}}(y) \cdot u(y)^2 \geq (b-a+1)N(a) + \sum_{d(x,y)=a} d_{\emph{\tiny out}}(y) \cdot u(y)^2.$$
If $G$ is locally contractive in the same region, then
$$  \sum_{d(x,y)=b+1} d_{\emph{\tiny in}}(y) \cdot u(y)^2 \leq (b-a+1)N(b) + \sum_{d(x,y)=a} d_{\emph{\tiny out}}(y) \cdot u(y)^2.$$
\end{corollary}

One way of phrasing the corollary is that in locally expansive regions, the size of $N(a)$ can be seen as a lower bound on the guaranteed $\ell^2$-growth of the harmonic function. In locally contractive regions, the size of $N(b)$ provides an upper bound on how much growth can have happened in that region. We note that, in contrast to Euclidean space, these estimates are weaker insofar as they provide only additive (as opposed to multiplicative) control; this reflects the fact that harmonic functions on graphs are more versatile than they are in Euclidean space (see \S 2.4 for an example).

\section{Proofs}

\subsection{Proof of the Theorem}

\begin{proof} For simplicity of exposition, we give the proof when the graph is simple, which means no self-loops and the edges are unweighted. The argument generalizes easily once the notion of in-degree and out-degree are refined to account for weights, see \S 2.5.
We fix $x \in V$ and start by showing that $N(0) \geq 0$. Note that
$$ N(0) = \left(\sum_{(x,y) \in E} u(y)^2\right) - \deg(x) \cdot u(x)^2.$$
Since the function is harmonic,
$$ u(x) = \frac{1}{\deg(x)} \sum_{(x,y) \in E} u(y),$$
which simplifies the problem to showing that
$$ \sum_{(x,y) \in E} u(y)^2 \geq \frac{1}{\deg(x)} \left( \sum_{(x,y) \in E} u(y) \right)^2.$$
This follows from the Cauchy-Schwarz inequality. It remains to prove the monotonicity.
We will work with the decomposition
$$ V = \bigcup_{k=0}^{\infty} V_k \qquad \mbox{where} \qquad V_k =\left\{w \in V: d(x,w) = k \right\}.$$
We will also introduce an abbreviation for edges that run between $V_k$ and $V_{k+1}$,
$$  E_k =\left\{e \in E: e~\mbox{runs between}~V_k~\mbox{and}~V_{k+1} \right\}.$$
For each vertex $y \in V_{k+1}$ there exist $y_1, \dots, y_{\ell} \subset V_k$ that are connected to $y$ by an edge. We have $\ell \geq 1$ because there exists a shortest path from $y$ to $x$. Then
\begin{align*}
    d_{\mbox{\tiny in}}(y) \cdot u(y)^2 &= \sum_{i=1}^{\ell} (u(y_i) + (u(y) - u(y_i)))^2  \\
   &\geq \sum_{i=1}^{\ell} u(y_i)^2 + 2 u(y_i) (u(y) - u(y_i)).
\end{align*}
We sum this inequality over all vertices $y \in V_{k+1}$ and obtain
\begin{align*}
    \sum_{y \in V_{k+1}} d_{\mbox{\tiny in}}(y) \cdot u(y)^2 &\geq \sum_{z \in V_k} d_{\mbox{\tiny out}}(z) u(z)^2 + \sum_{e = (z,y) \atop e \in E_k} 2 u(z) (u(y) - u(z)).
\end{align*}
We rewrite the second sum as
\begin{align*}
   \sum_{e = (z,y) \atop e \in E_k} 2 u(z) (u(y) - u(z)) &= \sum_{z \in V_k} \sum_{e = (z,y) \atop y \in V_{k+1}} 2 u(z) (u(y) - u(z)) \\
   &= \sum_{z \in V_k} 2 u(z) \sum_{e = (z,y) \atop y \in V_{k+1}}  (u(y) - u(z)).
\end{align*}
At this point, we use that $u$ is harmonic. For each $z \in V_k$, we have
$$ \sum_{(z,w) \in E} (u(w) - u(z)) = 0.$$
All the neighbors of $z \in V_k$ are either in $V_{k-1}$, in $V_k$ or in $V_{k+1}$ and thus
$$ \sum_{(z,w) \in E \atop w \in V_{k+1}}(u(w) - u(z)) = \sum_{(z,w) \in E \atop w \in V_{k-1}}(u(z) - u(w)) + \sum_{(z,w) \in E \atop w \in V_{k}}(u(z) - u(w)).$$
We first argue that, when summing over all $z \in V_k$, the second sum on the right-hand side can be simplified since we sum over each edge twice and
\begin{align*}
    \sum_{z \in V_k}  \sum_{(z,w) \in E \atop w \in V_{k}}  u(z) (u(z) - u(w))  &= \sum_{e =(a,b) \in E \atop a \in V_k,b  \in V_k} u(a) (u(a) - u(b)) +  u(b) (u(b) - u(a))\\
    &= \sum_{e =(a,b) \in E \atop a \in V_k,b  \in V_k}   (u(a) - u(b))^2 \geq 0.
\end{align*}
This implies
\begin{align*}
  \sum_{z \in V_k} 2 u(z) \sum_{e = (z,y) \atop y \in V_{k+1}}  (u(y) - u(z)) &\geq  
 \sum_{z \in V_k} 2 u(z) \sum_{e = (z,w) \atop w \in V_{k-1}}  (u(z) - u(w)) \\
 &=\sum_{z \in V_k} \sum_{e = (z,w) \atop w \in V_{k-1}} 2 u(z)  (u(z) - u(w)).
\end{align*}
Using
\begin{align*}
 2u(z)^2 - 2 u(z) u(w) \geq u(z)^2 - u(w)^2,
 \end{align*}
we have
$$\sum_{z \in V_k} \sum_{e = (z,w) \atop w \in V_{k-1}} 2 u(z)  (u(z) - u(w)) \geq \sum_{z \in V_k} \sum_{e = (z,w) \atop w \in V_{k-1}} u(z)^2 - u(w)^2.$$
This sum, in turn, is merely the sum over all edges $E_{k-1}$ and can be rewritten as
$$\sum_{z \in V_k} \sum_{e = (z,w) \atop w \in V_{k-1}} u(z)^2 - u(w)^2 = \sum_{z \in V_k} d_{\mbox{\tiny in}}(z) u(z)^2 -  \sum_{w \in V_{k-1}} d_{\mbox{\tiny out}}(w) u(w)^2.$$
Altogether, we see that
$$  \sum_{y \in V_{k+1}} d_{\mbox{\tiny in}}(y) u(y)^2 -  \sum_{z \in V_{k}} d_{\mbox{\tiny out}}(z) u(z)^2 \geq  \sum_{z \in V_k} d_{\mbox{\tiny in}}(z) u(z)^2 -  \sum_{w \in V_{k-1}} d_{\mbox{\tiny out}}(w) u(w)^2.$$
This is the desired statement.
 \end{proof}

\subsection{Proof of Corollary 1}
\begin{proof} The argument follows from telescoping the monotonicity formula.
We have
\begin{align*}
    \sum_{k=a}^{b} N(k) =  \sum_{k=a}^{b} \left[\sum_{d(x,y)=k+1} d_{\mbox{\tiny in}}(y) \cdot u(y)^2  - \sum_{d(x,y)=k} d_{\mbox{\tiny out}}(y) \cdot u(y)^2\right]. 
\end{align*}
If the graph is locally expansive, meaning 
$ d_{\mbox{\tiny out}}(y) \geq d_{\mbox{\tiny in}}(y)$,
then 
$$(b-a+1) N(a) \leq  \sum_{k=a}^{b} N(k) \leq \sum_{d(x,y)=b+1} d_{\mbox{\tiny in}}(y) \cdot u(y)^2 - \sum_{d(x,y)=a} d_{\mbox{\tiny out}}(y) \cdot u(y)^2.$$
Likewise, if the graph is locally contractive, meaning $ d_{\mbox{\tiny out}}(y) \leq d_{\mbox{\tiny in}}(y)$, then
$$(b-a+1) N(b) \geq  \sum_{k=a}^{b} N(k) \geq \sum_{d(x,y)=b+1} d_{\mbox{\tiny in}}(y) \cdot u(y)^2 - \sum_{d(x,y)=a} d_{\mbox{\tiny out}}(y) \cdot u(y)^2.$$
\end{proof}

\subsection{An example}
\begin{figure}[h!]
\begin{tikzpicture}
\filldraw (0,0) circle (0.08cm);
\filldraw (0,1) circle (0.08cm);
\draw [thick] (0,0) -- (0, 1);
\node at (0.3, 1) {0};
\node at (0, -0.3) {0};
\node at (1, -0.3) {1};
\node at (-1, -0.3) {-1};
\filldraw (1,0) circle (0.08cm);
\filldraw (2,0.5) circle (0.08cm);
\filldraw (2,-0.5) circle (0.08cm);
\filldraw (-1,0) circle (0.08cm);
\filldraw (-2,0.5) circle (0.08cm);
\filldraw (-2,-0.5) circle (0.08cm);
\draw [thick]  (-2, 0.5) -- (-1,0) -- (1,0) -- (2, 0.5);
\draw [thick]  (-2, -0.5) -- (-1,0);
\draw [thick] (1,0) -- (2, -0.5);
\draw [dashed] (2, 0.5) -- (3, 1);
\draw [dashed] (2, 0.5) -- (3, 0.5);
\draw [dashed] (2, -0.5) -- (3, -1);
\draw [dashed] (2, -0.5) -- (3, -0.5);
\draw [dashed] (-2, 0.5) -- (-3, 1);
\draw [dashed] (-2, 0.5) -- (-3, 0.5);
\draw [dashed] (-2, -0.5) -- (-3, -1);
\draw [dashed] (-2, -0.5) -- (-3, -0.5);
\draw [dashed] (0,1) -- (-0.5, 1.5);
\draw [dashed] (0,1) -- (0.5, 1.5);
\node at (0.8, 1.5) {0};
\node at (-0.8, 1.5) {0};
\end{tikzpicture}
\caption{Sketch of the function: a symmetric function around the root at 0, the third branch (going up) is set to always be 0.}
\label{fig:4}
\end{figure}
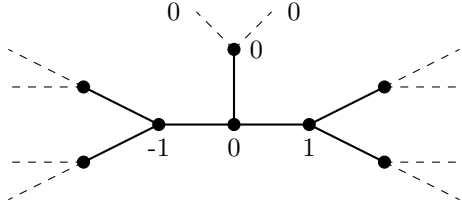

The geometry of infinite graphs can be quite different from that of the Euclidean space. Liouville's theorem can fail and bounded harmonic functions can exist, which leads to several other interesting consequences. The purpose of this section is to construct an explicit example of a harmonic function to illustrate our Theorem. We will work on the infinite 3-regular tree. 
Note that, except in the root, $d_{\mbox{\tiny in}}(v) = 1$ while $d_{\mbox{\tiny out}}(v) = 2$.
The function is sketched in Fig. \ref{fig:4}. We fix a vertex to be the root and set the function to be 0 there and in an entire branch leading away from it. We moreover set the harmonic function to be a function that is `symmetric' around the middle and merely a function of the distance to the origin. This leads to the sequence
$$ a_0 = 0, a_1 = 1 \qquad \mbox{and} \qquad a_{k+1} = \frac{3a_k - a_{k-1}}{2}.$$
One easily sees that, for $n \geq 1$
$$ a_n = 2 - \frac{1}{2^{n-1}}.$$

On the `left' branch, the values are simply $-a_n$. This leads to a non-constant harmonic function with values between $-2$ and $2$. We can now illustrate the discrete Almgren Monotonicity Formula for this function.
It states that
$$N(k) =  \sum_{d(x,y)=k+1} u(y)^2  - 2 \sum_{d(x,y)=k} u(y)^2 $$
is non-negative and that $N(k)$ is monotonically increasing in $k$. There are $2^{k}$ points at distance $k$ (counting both branches) and thus
$$ \sum_{d(x,y)=k} u(y)^2 = 2^k \left( 2 - \frac{1}{2^{k-1}} \right)^2.$$
A short computation shows that 
$$ N(k) =  2^{k+1} \left( 2 - \frac{1}{2^{k}} \right)^2 - 2^{k+1} \left( 2 - \frac{1}{2^{k-1}} \right)^2 = 8- \frac{3}{2^{k-1}}.$$
$N(k)$ non-negative and monotonically increasing and converges to 8.

\subsection{Weighted Graphs}
We briefly comment on the setting of weighted graphs $G=(V,E,w)$ where every edge $e \in E$ is additionally equipped with a weight $w_{e} > 0$. Harmonic functions are now functions $u:V \rightarrow \mathbb{R}$ satisfying
$$ \sum_{(x,y) \in E} w_{xy} (u(y) - u(x)) = 0.$$
We note that, in principle, we could even allow for self-loops $(x,x) \in E$ since they automatically disappear when considering harmonic functions. 
In the weighted setting, we can adapt the notion of in-degree and out-degree for $v\in V_k$  as
$$ d_{\mbox{\tiny in}}(v) = \sum_{e = (x,v) \in E \atop x \in V_{k-1}} w_{xv} \qquad \mbox{and} \qquad d_{\mbox{\tiny out}}(v) =  \sum_{e = (v,x) \in E \atop x \in V_{k+1}} w_{vx}.$$
Note that if the edges all have weight 1, then these definitions reduce themselves to the definitions already used above. Having these definitions in place, we will conclude just as above that
$$N(k) =  \sum_{d(x,y)=k+1} d_{\mbox{\tiny in}}(y) \cdot u(y)^2  - \sum_{d(x,y)=k} d_{\mbox{\tiny out}}(y) \cdot u(y)^2 $$
is non-negative and satisfies $N(k+1) \geq N(k)$. Showing that $N(1) \geq N(0)$ then amounts to showing that
\begin{align*}
     \sum_{d(x,y) = 1} w_{xy} u(y)^2 \geq \left(\sum_{d(x,y) = 1} w_{xy} \right) u(x)^2 =    \frac{1}{d_{\mbox{\tiny out}}(x)} \left(\sum_{d(x,y) = 1} w_{xy} u(y) \right)^2 
\end{align*}
which follows from applying Cauchy-Schwarz to
$$ \left(\sum_{d(x,y) = 1} w_{xy} u(y) \right)^2  = \left(\sum_{d(x,y) = 1} \left(\sqrt{w_{xy}} u(y)\right) \cdot \sqrt{w_{xy}} \right)^2.$$
As for the remainder of the argument, we argue in parallel and obtain
\begin{align*}
    \sum_{y \in V_{k+1}} d_{\mbox{\tiny in}}(y) \cdot u(y)^2 &\geq \sum_{z \in V_k} d_{\mbox{\tiny out}}(z) u(z)^2 + \sum_{e = (z,y) \atop e \in E_k} 2 w_{e} u(z) (u(y) - u(z)).
\end{align*}
The weight $w_e$ can then be carried through the rest of the argument without having any further impact and the result follows.


\begin{thebibliography}{10}
\bibitem{almgren} F. Almgren, Dirichlet's problem for multiple valued functions and the regularity of mass minimizing integral currents. Minimal submanifolds and geodesics, Proc. Japan–United States Sem., Tokyo, 1977, North-Holland, Amsterdam, New York (1979), pp. 1--6

\bibitem{almgren2} F. Almgren, Almgren’s big regularity paper, World Scientific Monograph Series in Mathematics, 1, World Scientific Publishing Co. Inc., River Edge, NJ, 2000.

\bibitem{ACS} I. Athanasopoulos, L. A. Caffarelli, and S. Salsa, The structure of the free boundary for
lower dimensional obstacle problems. Amer. J. Math. 130 (2008), no. 2, p. 485--498.

\bibitem{graph0} I. Benjamini and O. Schramm, Harmonic functions on planar and almost planar graphs and manifolds, via circle packings. Inventiones mathematicae, 126 (1996), 565-587.

\bibitem{graph1} M. Guadie and E. Malinnikova, On three balls theorem for discrete harmonic functions. Computational Methods and Function Theory, 14 (2014), p. 721--734.

\bibitem{graph2} H. Heilbronn, On discrete harmonic functions,  Mathematical Proceedings of the Cambridge Philosophical Society 45 (1949), p. 194-206

\bibitem{graph3} L. Lov\'asz,  Discrete analytic functions: an exposition. Surveys in differential geometry 9 (2004), p. 241--273.


\bibitem{PSU} A. Petrosyan, H. Shahgholian and N. Uraltseva, Regularity of Free Boundaries in Obstacle-Type Problems, Graduate Studies in Mathematics, Amer. Math. Soc. 136 (2012).

\bibitem{sp} L. Spolaor, Monotonicity Formulas in the Calculus of Variation. Notices of the American Mathematical Society 69 (2022), p. 1731-1737.

\end{thebibliography}
\end{document}